
\input amstex
\documentstyle{amsppt}
\catcode34=12 
\magnification \magstep1 \TagsOnRight \pageheight{24truecm}
\pagewidth{17truecm} \pageno1 \nopagenumbers

\def\myfootline{\hss\tenrm\folio}
\headline={}
\footline={\ifnum\pageno>1\myfootline\fi} \NoRunningHeads

\expandafter\redefine\csname logo\string @\endcsname{}

\textfont2=\tensy \scriptfont2=\sevensy \scriptscriptfont2=\fivesy
\def\cal{\fam2}
\def\eqalign{\null\vbox\bgroup\advance\baselineskip3pt\halign
\bgroup\hfil${}##{}$&&${}##{}$\hfil\crcr}
\def\endeqalign{\crcr\egroup\egroup\,}

\vskip1cm \topmatter

\newcount\chapno    
\newcount\claimno   
\newcount\stateno   
 \chapno=0
 \stateno=0
 \claimno=0
\define\chap{\global\claimno=0\global\stateno=0\global\advance\chapno
by 1
\ \the\chapno.\enskip}

\define\Clabel{{\global\advance\claimno by
1}\number\chapno.\number\claimno}
\define\clabel#1{\global\advance\claimno by 1\expandafter\xdef
 \csname
Cl#1\endcsname{\the\chapno.\the\claimno}\the\chapno.\the\claimno}

\define\Elabel{\tag\global\advance\stateno by
1\the\chapno.\the\stateno}
\define\elabel#1{\tag\global\advance\stateno by 1\expandafter\xdef
 \csname St#1\endcsname{\thetag
       {\the\chapno.\the\stateno}}\the\chapno.\the\stateno}

 \define\cref#1{\expandafter\csname Cl#1\endcsname}
 \define\eref#1{\expandafter\csname St#1\endcsname}

\define\mps#1#2#3{#1\colon #2\to #3}
\define\Mps#1#2#3{#1\colon #2\hookrightarrow #3}
\define\pam#1#2#3#4{#1\hookleftarrow #2@>#3>>#4}

\def\>{\prec}
\def\<{\succ}

\def\dim{\op{dim}}

\def\E{{\cal E}}

\def\L{{\cal L}}

\def\O{{\cal O}}
\def\P{{\cal P}}

\def\U{{\cal U}}
\def\V{{\cal V}}

\define\op#1{\operatorname {#1}}
\def\dist#1#2{\op{dist}(#1,#2)}

\def\dim{\op{dim}}

\def\Int{\op{Int}}

\def\mean{\rightleftharpoons}

\def\tilde{\widetilde}
\def\phi{\varphi}
\def\matop#1#2#3{\mathop{#1}\limits_{#2}^{#3}}

\def\Id{\op{Id}}
\def\Bd{\op{Bd}}

\def\ANE{\op{ANE}}
\def\AE{\op{AE}}

\def\Eucl{\op{Eucl}}

\topmatter \abstract We prove that for every $n>2$, the
Banach-Mazur compactum $Q(n)$ is the compactification of a Hilbert
cube manifold by the Euclidean point. For $n=2$ this result was
proved earlier.
\endabstract
\thanks
The first author was supported in part by the NSERC grant No.
OGP005616. The second author was supported in part by the RFFI
grant No. 00-01-00289. The third author was supported in part by
the MESSRS grant No. 0101-509.
\endthanks
\author S.M. Ageev, S.A. Bogatyi and D. Repov\v{s}
\endauthor
\title
The Banach-Mazur compactum is the Alexandroff compactification of
a Hilbert cube manifold
\endtitle
\rightheadtext{On Banach-Mazur compacta}
\address
Department of Mathematics and Statistics, University of
Saskatchewan, Saskatoon, Canada  SS  S7N 5A5. email:
ageev\@snoopy.usask.ca\endaddress
\address Faculty of Mechanics and Mathematics,
Moscow State University, Vorobyovy Gory, \break Moscow, Russia
119992. email: bogatyi\@mech.math.msu.su \endaddress
\address Institute of Mathematics, Physics and Mechanics,
University of Ljubljana, P. O. Box 2964, Ljubljana, Slovenia 1001.
email: dusan.repovs\@uni-lj.si\endaddress
\endtopmatter

\document
\head\chap Introduction
\endhead
In this paper we continue the study of topological properties of
the Banach-Mazur compacta $Q(n)$, i.e. the sets of all isometry
classes of $n$-dimensional Banach spaces, topologized by the
Banach-Mazur metric. Recently, a substantial progress has been
made concerning these spaces. It was proved in \cite{ABF} that
$Q(n)\in \op{AE}$ for all $n\ge 2$. The long standing problem
about topological equivalence of $Q(n)$ and the Hilbert cube
$I^{\infty}$ was finally solved in the negative in \cite{AB} (for
details see \cite{AB1} and \cite{R}): \proclaim{Theorem
\clabel{th}} $Q(2)$  and $I^{\infty}$ are not homeomorphic.
\endproclaim
The key idea of the proof of Theorem \cref{th} is to show that
$Q(2)\setminus\{\Eucl\}$ fails to be homotopically trivial where
$\{\Eucl\}\in Q(n)$ is the Euclidean point which corresponds to
the isometry class of the standard $n$-dimensional Euclidean
space. In turn, this is a corollary of non-triviality of the
four-dimensional cohomological group with rational coefficients
$\op{H}^4(Q(2)\setminus\{\Eucl\},\Bbb Q)$. An intimate connection
between the Banach-Mazur compacta and the Smith theory of periodic
homeomorphisms \cite{Br} was first revealed in \cite{AB1; p.7} and
\cite{AB}. The same argument was used in \cite{An1; Corollary 6,
p.224}.

The investigations of $Q(2)$  were continued in \cite{AR1} -- it
was proved that $Q(2)$ is the one-point compactification of a
Hilbert cube manifold. This implies, with the help of Theorem
\cref{th}, the non-homogeneity of $Q(2)$. A natural problem about
the structure of $Q(n),n>2$, was reduced to a plausible question
from convex geometry. Here we give the affirmative (and complete)
answer to this problem: \proclaim{Theorem \clabel{Th1+0}} $Q_{\cal
E}(n)=Q(n)\setminus \{\Eucl\}$ is a $I^{\infty}$-manifold.
\endproclaim

\head\chap Preliminary facts and results
\endhead
The space $Q(n)$ of all isometry classes of $n$-dimensional Banach
spaces is called the {\it Banach-Mazur compactum} (see \cite{AR1}
and \cite{R} for detailed exposition of topology of Banach-Mazur
compacta). This compactum admits a representation as a
decomposition of the space $C(n)$ of all compact convex bodies in
$\Bbb R^n$ which are symmetric $\op{rel}\ 0$.

If one measures the distance between subsets of $\Bbb R^n$ by the
Hausdorff metric $\rho_H$ and defines the linear combination
$\sum_{i=0}^n\lambda_iA_i$ by means of the Minkowski operation,
then $(C(n),\rho_H)$ becomes a locally compact convex space.
Moreover, $C(n)$ can be equipped by an action of the general
linear group $$\op{GL}(n)\times C(n)\to C(n),T\cdot V=T(V),$$
where $T:\Bbb R^n\to\Bbb R^n\in GL(n)$ and $V\in C(n)$, which
agrees with the convex structure on $C(n)$. It is well-known that
the orbit space $C(n)/GL(n)$ is naturally homeomorphic to the
Banach-Mazur compactum.

Let $G$ be a compact Lie group. An {\it action} of $G$ on a space
$X$ is a homomorphism $T:G\to\op{Aut} X$ of the group $G$ into the
group $\op{Aut} X$ of all autohomeomorphisms of $X$, such that the
map $G\times X\to X$, given by $(g,x)\mapsto T(g)(x)=g\cdot x$, is
continuous. A space $X$ with a fixed action of $G$ is called a
$G$-{\it space}.

For any point $x\in X$, the {\it isotropy subgroup} of $x$, or the
{\it stabilizer} of $x$, is defined as $ G_x=\{g\in G\mid g\cdot
x=x\}$ and the {\it orbit} of $x$ as $ G(x)=\{g\cdot x\mid g\in
G\}$. The space of all orbits is denoted by $X/G$ and the natural
map $\pi:X\to X/G$, given by $\pi(x)=G(x)$, is called the {\it
orbit projection}. The orbit space $X/G$ is equipped with the
quotient topology, induced by $\pi$. For more details see
\cite{Br}.

A space $X$ is called an {\it absolute neighborhood extensor},
$X\in \ANE$, if every map $\phi:A\to X$, defined on a closed
subset $A\subset Z$ of a metric space $Z$, and called a {\it
partial map}, can be extended over some neighborhood $U\subset Z$
of $A$, $\tilde{\phi}:U\to X, \tilde{\phi}\restriction_A=\phi$. If
we can always take $U=Z$ then $X$ is called an {\it absolute
extensor}, $X\in \AE$. 

We note that in the case when $X$ is a
metric space, the concepts of the absolute (neighborhood) retract
and the absolute (neighborhood) extensor coincide \cite{Hu},
\cite{M}. By the Toru\'nczyk Characterization Theorem \cite{T}, a locally
compact space $X\in\ANE$ is a $I^\infty$-manifold if and only if
$X$ admits arbitrary small maps $f_i:X\to X$, $i\in\{1,2\}$, with
$\op{Im} f_1\cap \op{Im} f_2=\emptyset$.

 Let $(X,d)$ be a metric space of diameter 1.
The following formula for a metric on the cone $\op{Con} X$ was
quoted in \cite{AR1}:
$$\rho((x,t),(x',t'))=\sqrt{t^2+(t')^2-2tt'\cos\gamma},
\hbox{where }\cos\gamma=(2-d^2(x,x'))/2.$$ However, the correct
formula,
which is only slightly different, should have been the following
well-known one (see e.g. \cite{BBI; p.91}):
$$\rho((x,t),(x',t'))=\sqrt{t^2+(t')^2-2tt'\cos\gamma},
\hbox{where } \gamma=d(x,x').$$ The authors acknowledge S. A.
Antonyan for kindly pointing out to us this {\it lapsus} - more on
this formula can be found in \cite{An2}.

Next, we introduce a partial order among compact Lie groups. We
set $L<H$, where $L$ and $H$ are compact Lie groups, if $L$ is
isomorphic to a proper subgroup of $H$. Clearly,
\item\item{$(\alpha)$}
If $L<H$, then either $\dim L<\dim H$, or $\dim L=\dim H$ and
${\cal C}_L<{\cal C}_H$ where ${\cal C}_H$ is the number of
components of linear connectivity of $H$.

  If $L$ is a closed subgroup of the compact
Lie group $H$ and $\dim L=\dim H$, then $L$ is an open subgroup.
 The
 following stronger fact can easily be derived
 from this:
 \item\item{$(\beta)$}
 Let $L$ be a closed subgroup of the compact
Lie group $H$. If $\dim L=\dim H$ and ${\cal C}_L={\cal C}_H$,
then $L=H$.

 The verification of
 the
 following property of the
 order
 introduced above
 is quite
easy, and so we leave it to the reader.
 \item\item{$(\gamma)$}
 There does not exist
 any
countable sequence of compact Lie groups $\{H_i\}$ such that
$H_1>H_2>H_3>\dots>H_n>\dots$.

 Clearly, the pair $\op{ind}H=(\dim H,{\cal C}_H)$
 belongs to
$\Bbb N\times \Bbb N$.  Let $\Bbb N\times \Bbb N$ be endowed with
the lexicographical order. It follows from $(\alpha)$ that the map
$H\mapsto\op{ind}H$ preserves the order.

We can derive the following principle from $(\gamma)$  permitting
us to prove by induction for compact Lie groups:
 \proclaim{Proposition \clabel{t5.9++}}
 Let $\P(H)$ be a property
which depends on a compact Lie group $H$. Suppose that
  \item\item{$(\delta)$} $\P(H)$ is true for the trivial group
$H=\{e\}$; and
 \item\item{$(\varepsilon)$} $\P(H)$ is true if
$\P(L)$ is true for every $L<H$.

\noindent Then $\P(H)$ is true for every group $H$.
\endproclaim
 For instance,
 we examined
 in \cite{AR2}  the following property
satisfying $(\delta)$ and $(\varepsilon)$: "$\P(H)$ is valid if
and only if for every metric $H$-space $X\in H$-$\ANE$, the orbit
space $X/H$ is $\ANE$".

 It is well-known that (see \cite{J}) for every convex body $V\in
C(n)$, there exists a unique ellipsoid $E_V\in C(n)$ (called the
{\it L\"owner ellipsoid}), which contains $V$ and has the minimal
Euclidean volume. The minimality of $\op{vol} E_V$ implies the
$GL(n)$-invariance of $E_V$, $E_{T\cdot V}=T\cdot E_V$, for every
$T\in\op{GL}(n)$. A continuous dependence $E_V$ of $V$ with
respect to Hausdorff metric was proved in \cite{ABF}. Therefore
${\cal L}:C(n)\to\goth E$, ${\cal L}(V)=E_V$, is a
$\op{GL}(n)$-retraction of $C(n)$ onto the ellipsoid orbit $\goth
E=\op{GL}(n)\cdot B^n$ where $B^n$ is the unit ball (${\cal L}$ is
said to be the {\it L\"owner retraction}).

In the sequel, $O(n)$ will denote the orthogonal group of $\Bbb
R^n$. Let $L(n)={\cal L}^{-1}(B^n)$ be the $O(n)$-slice which is
an $O(n)$-space. In other words, $L(n)$ consists of all bodies
$V\in C(n)$ whose minimal L\"owner ellipsoid coincides with $B^n$.
The orbit space $Q(n)=C(n)/GL(n)$ is homeomorphic to $L(n)/O(n)$.
Since, by \cite{ABF},  $Q(n)\in \AE$, it follows that
$L(n)/O(n)\in \AE$, and therefore $Q_{\cal E}\mean
Q(n)\setminus\{\Eucl\} =L_{\cal E}(n)/O(n)\in \ANE$ where $L_{\cal
E}\mean L(n)\setminus\{B^n\}$. Hence, Theorem \cref{Th1+0} is
reduced to the following assertion:
 \proclaim{Theorem
\clabel{Th1+00}} $L_{\cal E}(n)/O(n)$ is an $I^{\infty}$-manifold.
\endproclaim
It is well known \cite{Ab} that there exists an $O(n)$-retraction
${\cal R}:C(n)\to L(n) $ which maps $C_{\cal E}(n)$ into $L_{\cal
E}(n)$. However, we need the following folklore result which gives
more information and follows from geometric considerations (see
e.g. \cite{AR1}):
 \proclaim{Proposition \clabel{t5.9}}
 There
exists a continuous $O(n)$-retraction $\goth R:C(n)\to L(n)$ such
that for every $V\in C(n)$, $\goth R(V)$ and $V$ are affinely
equivalent.
\endproclaim
\demo{Hint} Let $T$ be an element of $\op{GL}(n)$ such that
$T^{-1}\cdot B^n=\L(V)$. By \cite{L}, $T$ can be represented as
$T_2\circ T_1$ where $T_2\in O(n)$ and $T_1$ is self-adjoint. We
set $\goth R(V)=T_1(V)$.  We leave the verification of the
required properties to the reader.
\enddemo

\head\chap Proof of Theorem \cref{Th1+00}
\endhead
The proof of Theorem \cref{Th1+00} (and therefore also Theorem
\cref{Th1+0}) is easily reduced - by invoking the Toru\'nczyk
Characterization Theorem - to the following fact (see \cite{AR1}):
\proclaim{Theorem \clabel{pr5.3}} For every $\delta>0$, there
exist $O(n)$-maps $f_i:L_{\cal E}(n)\to L_{\cal E}(n)$,
$i\in\{1,2\}$, such that $\dist{f_i}{\Id_{L_{\cal E}(n)}}<\delta$
and $\op{Im} f_1\cap \op{Im} f_2=\emptyset$.
\endproclaim
In turn, Theorem \cref{pr5.3} evidently reduces to the following
two theorems. We first recall some necessary notions. A point $a$
of a convex set $V\subset\Bbb R^n$ is called {\it extreme} if
$V\setminus\{a\}$ is convex. It is well-known that the set
$\op{Extr}(V)$ of all extreme points of $V$ lies in the relative
boundary $\op{rbd}V$, and $V$ coincides with the convex hull
$\op{Conv}(\op{Extr}(V))$ of $\op{Extr}(V)$. If
$\op{Extr}(V)=\op{rbd}(V)$, $V$ is called {\it elliptically
convex}, otherwise $V$ is called {\it non-elliptically convex}.
 \proclaim{Theorem \clabel{Th1}}
 There
exists an $O(n)$-homotopy $H:L(n)\times [0,1]\to L(n)$ such that:
\item\item{\rm(a)} $H_0=\Id$;
\item\item{\rm(b)} If $V\in L(n)$ and $t\in[0,1]$, then
$H_t(V)=B^n$ if and only if $V=B^n$; and
\item\item{\rm(c)}  $H_t(V)$ is elliptically convex
for each $V\in L(n)$ and each $t>0$.
\endproclaim
\proclaim{Theorem \clabel{Th2}} There exists an $O(n)$-homotopy
$F:L(n)\times [0,1]\to L(n)$ such that:
\item\item{\rm(d)} $F_0=\Id$; and
\item\item{\rm(e)}  $F_t(V)$ is non-elliptically convex
for each $V\in L_{\cal E}(n)$ and each $t>0$.
\endproclaim
{\it Proof of Theorem \cref{Th1}.} Let $V\in L(n)$. Note that a
convex body $V$ is elliptically convex if and only if every
supporting hyperplane of $V$ and $V$ are intersected at one point
\cite{L}. This criterion permits to match the notion of elliptical
convexity with the Minkowski linear combination.
 \proclaim{Lemma \clabel{l}}
 Let
$V,V_i\in C(n)$ and $V=\matop{\sum}{i=1}{p}\lambda_i\cdot V_i$,
where all $\lambda_i>0$. Then $V$ is elliptically convex if and
only if $V_i$ is elliptically convex for each $i$.
\endproclaim
\demo{Proof} Let $\Pi$ and $\Pi_i$ be parallel supporting closed
hyperplane of $V$ and $V_i$, respectively. Let also $A=V\cap\Pi$
and $A_i=V_i\cap\Pi_i$. It is evident that
$A=\matop{\sum}{i=1}{p}\lambda_i\cdot A_i$, and
$\matop{\sum}{i=1}{p}\lambda_i\cdot x_i\in A,x_i\in V_i$ if and
only if $x_i\in A_i$ for each $i$. Hence, $A$ consists of one
point if and only if each $A_i$ consists of one point. Now apply
the criterion of elliptical convexity mentioned above. $\
\blacksquare$
\enddemo
Let $\Psi:\Bbb R^n\times [0,1]\rightarrow\Bbb R^n$ be defined by
$$\Psi(x,t)=\Psi_t(x)\mean x/(1+t\cdot \|x\|)\in\Bbb R^n.$$
Clearly, for every $t\in [0,1]$
\item\item{\rm(1)} $\Psi_t$ is a
continuous $O(n)$-embedding; and
\item\item{\rm(2)}
If $\Psi_{t}(V)$ is affinely equivalent to $B^n$ for some $V\in
L(n)$, then $V=B^n$.

\proclaim{Lemma \clabel{Th1+1}}(see \cite{BP; p.95}) For any
$t\in(0,1]$ and for any $V\in C(n)$, $\Psi(V,t)$ is an
elliptically convex body.
\endproclaim
 Clearly, the $O(n)$-map
$$ \mps{\Psi}{L(n)\times [0,1]}{C(n)},(V,t)\in L(n)\times
[0,1]\matop{\mapsto}{}{} \Psi(V,t)\in C(n),$$ is continuous. As an
easy corollary of Lemma \cref{Th1+1} and $(1)-(2)$ we conclude
that the $O(n)$-homotopy $\mps{H\mean \goth R\circ\Psi}{L(n)\times
[0,1]}{L(n)}$, where $\goth R$ is taken from Proposition
\cref{t5.9}, satisfies Theorem \cref{Th1}. Details are left to the
reader.

{\it Proof of Theorem \cref{Th2}.} First, observe that no
additional extreme points appear if we take a convex hull of the
union of convex sets.
\item\item{\rm(3)}
$\op{Extr}(\op{Conv}A)\subset A$, for every $A\subset\Bbb R^n$.

 Hence, if $L$ is a finite subset, then
$\op{Extr}\op{Conv}L$ is also finite, and therefore $\op{Conv}L$
is non-elliptically convex. The key argument in the proof of
Theorem \cref{Th2} is concerned with finding more non-elliptically
convex bodies in $C(n)$.
 \proclaim{Lemma \clabel{l5.8}}
 Let $V\in L_{\E}(n)$ and let
$H=O(n)_V$ be a stabilizer of $V$. Then for every finite subset
$L\subset\Bd V$ with $0\in\Int(\op{Conv}L)$,
$W\mean\op{Conv}(H\cdot L)\in C(n)$ is non-elliptically convex.
Moreover, $O(n)_W\supset H$.
\endproclaim
\demo{Proof} Since $\{\pm\Id\}\subset H$ and
$0\in\Int(\op{Conv}L)\subset\op{Conv}(H\cdot L)$, it follows that
$\op{Conv}(H\cdot L)\in C(n)$.

We note the following evident fact:
\item\item{$(4)$} If $A,B\in C(n)$ and $\Bd
A\subset \Bd B$, then $A=B$ (and therefore $\Bd A=\Bd B$).

Suppose, contrary to the assertion of the lemma, that $W$ is
elliptically convex, i.e. $\break\op{Extr}(W)=\Bd W$. Therefore
$$\Bd W=\op{Extr}(\op{Conv}(H\cdot L))\subset H\cdot L\subset \Bd
V.$$ Hence, by $(4)$,  $V=W$. Since $L$ is finite and $O(n)$ acts
orthogonally on $\Bbb R^n$, $H\cdot L$ is contained in the disjoint
union $\bigsqcup r_i\cdot S^{n-1}$ of finitely many concentric
spheres. In view of $\Bd W\subset H\cdot L$ and connectivity of
$\Bd W$, $\Bd W\subset r_{i_0}\cdot S^{n-1}$ for some $i_0$. By
$(4)$, we have $V=W=r_{i_0}\cdot B^{n}$ which contradicts with
$V\in L_{\E}(n)$.

 Finally, $O(n)_W$
contains $H$,
 in view of
the orthogonality of the action of $O(n)$ on $\Bbb R^n$. $\
\blacksquare$
\enddemo
Next, we find in arbitrary neighborhood of $V\in L_{\E}(n)$ a
non-elliptically convex body with more properties.
\proclaim{Proposition \clabel{l5.9}} Let $V\in L_{\E}(n)$ and let
$H=O(n)_V$ be a stabilizer of $V$. Then for every $\varepsilon>0$
there exists a finite set $L\subset\Bd V$ with
$0\in\Int(\op{Conv}L)$ such that:
\item\item{\rm(i)} $W=\op{Conv}(H\cdot L)\in C(n)$ is
non-elliptically convex;
\item\item{\rm(ii)} $V$ and $W$ have equal stabilizers; and
\item\item{\rm(iii)} $\rho_H(V,W)<\varepsilon$.
\endproclaim
\demo{Proof} We apply \cite{Br; 5.5} for the  $O(n)$-space $C(n)$.
Then there exists a $\theta>0$ such that $U\in
C(n),\rho_H(V,U)<\theta$ implies:
\item\item{\rm(5)} The conjugate subgroup to $O(n)_{U}$ is
contained in $H$, or,  equivalently, $O(n)_{U}\subset H'$, where
$H'$ is a subgroup conjugate to $H$.
 \proclaim{Lemma \clabel{l5.9+0}}
 If
$U\in C(n), \rho_H(V,U)<\theta$ and $O(n)_{U}\supseteq H$, then
$O(n)_{U}=H$.
\endproclaim
\demo{Proof} Since $H\subset O(n)_{U}$, we have, by $(5)$, that
$H\subset H'$. Since the Lie groups $H$ and $H'$ are isomorphic,
$\dim H=\dim H'$ and ${\cal C}(H)={\cal C}(H')$. The property
$(\beta)$ implies that $H=H'$. $\ \blacksquare$
\enddemo
 Clearly, there exists a finite subset $L\subset\Bd V$ such that
$0\in\Int(\op{Conv}L)$ and $\break\rho_H(V,\op{Conv}L)<\theta$. By
Lemma \cref{l5.8}, $W\mean\op{Conv}(H\cdot L)\in C(n)$ is
non-elliptically convex. Then $$V\supseteq\op{Conv}(H\cdot
L)=W\supseteq \op{Conv}L,$$ and therefore $\rho_H(V,W)<\theta$.
Since $$O(n)_W= O(n)_{\op{Conv}(H\cdot L)}\supset O(n)_{H\cdot
L}\supseteq H,$$ we have by Lemma \cref{l5.9+0}, that $O(n)_W=H$.
$\ \blacksquare$
\enddemo
Now, we apply Proposition \cref{l5.9} to construct an abundant
collection of equivariant retractions onto non-elliptically convex
orbits. We say that a family $\{ B_{\gamma}\}$ of $G$-subsets of
$X$ lying in $X\setminus C$ is {\it $G$-adjoint to the $G$-set
$C$}, provided that for every $x\in\ C$ and for every neighborhood
$\O(x)\subset X$, there exists a neighborhood $\O_1(x)\subset X$
such that $B_\gamma\subset G\cdot\O(x)$ as soon as $B_\gamma\cap
G\cdot \O_1(x)\not=\varnothing$.
 \proclaim{Lemma \clabel{Th2+1}}
There exist an open $O(n)$-cover $\omega=\{\U_\gamma\}$ of $Z\mean
L_{\cal E}(n)\times(0,1]$ and a family
$\Omega=\{r_\gamma:\U_\gamma\to P_\gamma\}$ of $O(n)$-maps such
that:
\item\item{\rm(f)} For each $\gamma$,
$P_\gamma=(O(n)\cdot Q_\gamma)\times\{t_\gamma\}$, where
$Q_\gamma\in C(n)$ is a non-elliptically convex body and
$t_\gamma\in(0,1]$;
\item\item{\rm(g)}
$\{\U_\gamma\}$ is an $O(n)$-adjoint cover to $$A\mean L(n)\times
[0,1]\setminus Z=\{B^n\}\times [0,1]\cup L(n)\times \{0\};$$ and
\item\item{\rm(h)}
For every $a\in A$ and every $\varepsilon>0$, there exists
$\delta>0$ such that
$\dist{r_{\gamma}}{\Id_{\U_\gamma}}<\varepsilon$ as soon as
$\U_{\gamma}$ is contained in the $\delta$-neighborhood (with
respect to Hausdorff metric) $\op{N}(\O(a);\delta)$ of the
$\op{O}(n)$-orbit $\O(a)\subset A$ of $a$ (or briefly,
$\dist{r_{\gamma_i}}{\Id}\rightarrow 0$, whenever
$\U_{\gamma_i}\rightarrow A$).
\endproclaim
\demo{Proof} Let $Q\in L_{\cal E}(n)$, $t\in (0,1]$, and let
$R=\{(g\cdot Q,t)|g\in O(n)\}$ be an orbit of $Z$. By the Palais
Slice theorem \cite{P}, there exists an $O(n)$-retraction $r'_R:
{\V}_R\to R$, $r'_R\restriction_{R}=\Id$, where ${\V}_R\subset Z$
is an invariant neighborhood of $R$. Here we can assume that:
\item\item{$(6)$} $\{\V_R\}$ is an $O(n)$-adjoint cover to
$A$; and
\item\item{$(7)$}
$\dist{r'_{R_i}}{\Id}\rightarrow 0$, whenever $\V_{R_i}\rightarrow
A$.

By Proposition \cref{l5.9}, for every orbit $R=(O(n)\cdot
Q)\times\{t\}\subset Z,$ we can fix an orbit $R'=(O(n)\cdot
Q')\times\{t\}$ such that $Q'\in C(n)$ is non-elliptically convex,
and there exists an $O(n)$-homeomorphism $s_R:R\rightarrow R'$
with $\dist{s_{R_i}}{\Id}\rightarrow 0$ whenever
$\V_{R_i}\rightarrow A$. The cover $\{\V_R\}$ and the family of
compositions $r_R=s_R\circ r'_R: {\V}_R\to R'$ are the desired
objects $\omega$ and $\Omega$. $\blacksquare$
\enddemo

We now complete the proof of Theorem 3.3. Let
$\{\lambda_\gamma:Z\to[0,1]\}$ be a continuous equivariant
partition of unity, subordinate to the cover
$\omega=\{\U_\gamma\}$. Let $\goth R$ be a retraction from
Proposition \cref{t5.9}. We define the desired $O(n)$-map
$F:L(n)\times [0,1]\to L(n)$ as follows:
 $$F(V,t)=\cases \goth R\circ
(\matop{\sum}{\gamma}{}\lambda_\gamma(V,t)\cdot r_{\gamma}(V,t)),\
\ \text{where}\ (V,t)\in Z,\ \text{and}\cr
 F(V,t)=V,\ \ \ \ \ \ \ \ \ \ \ \ \ \ \ \ \ \
 \text{where}\ (V,t)\in A. \endcases $$
   By Lemma \cref{l},
$\matop{\sum}{\gamma}{}\lambda_\gamma(V,t)\cdot
r_{\gamma}(V,t),(V,t)\in Z$, is non-elliptically convex. Since, by
Proposition \cref{t5.9}, $\goth R(W)$ and $W$ are affinely
equivalent, it follows that  $F(V,t),(V,t)\in Z$, is also
non-elliptically convex. The continuity of $F$ at a point
$(V,t)\in A$ follows from $(7)$. $\blacksquare$

\medskip\medskip\medskip\medskip

\enddocument